\documentclass[10pt,a4paper,reqno]{amsart} 
\usepackage{latexsym}
\usepackage{verbatim}
\usepackage{epsfig}
\usepackage{rotating}
\usepackage{amssymb}
\usepackage{afterpage}
\usepackage{color}
\usepackage{dsfont}

\usepackage{amssymb,amsfonts,amsmath,stmaryrd}

\newtheorem{Theorem}{Theorem}

\newtheorem{Proposition}[Theorem]{Proposition}
\newtheorem{Definition}[Theorem]{Definition}

\catcode`\@=11
\def\section{\@startsection{section}{1}%
 \z@{.7\linespacing\@plus\linespacing}{.5\linespacing}%
 {\normalfont\bfseries\scshape\centering}}

\def\subsection{\@startsection{subsection}{2}%
  \z@{.5\linespacing\@plus\linespacing}{.5\linespacing}%
  {\normalfont\bfseries\scshape}}

\def\subsubsection{\@startsection{subsubsection}{3}%
 \z@{.5\linespacing\@plus\linespacing}{-.5em}
  {\normalfont\bfseries\itshape}}
\catcode`\@=12

\newcommand{\ns}{\mathbb{N}}

\def\emm#1,{{\em #1}}

\newcommand{\la}{\lambda}

\begin{document}
\title
[A  property of  dominance on partitions]
{A    property of dominance  on partitions}

\author[M. Bousquet-M\'elou]{Mireille Bousquet-M\'elou}
\address{M. Bousquet-M\'elou: CNRS, LaBRI, Universit\'e Bordeaux 1, 
351 cours de la Lib\'eration, 33405 Talence, France}   
\email{mireille.bousquet@labri.fr}
\date{March 18, 2008}


\begin{abstract}
Given an integer  partition $\la=(\la_1, \ldots, \la_\ell)$ and an integer $k$,  denote by $\la^{(k)}$  the sequence of length $\ell$
obtained by reordering the values $|\la_i-k|$ in non-increasing
order.  If $\la$ dominates $\mu$ and has the same weight, then  $\la^{(k)}$
  dominates $\mu^{(k)}$.
\end{abstract}

\maketitle


A \emm partition, of an integer $n$ in $k$ parts is a non-increasing
sequence  $\la=(\la_1, \ldots,  \la_k)$ of positive integers that sum
to $n$. The integer $n$ is called the \emm weight, of $\la$,
denoted $||\la||$, 
while $k$ is the \emm length, of $\la$, denoted  $\ell(\la)$. Partitions 
are partially ordered by dominance: 
$$
\la\ge\mu \quad \hbox{ if} \quad \la_1+\cdots +\la_i \ge \mu_1+\cdots
+\mu_i 
\quad \hbox{ for all } i,
$$
where it is understood that $\la_i=0$ if $i>\ell(\la)$ (and similarly
for  $\mu$).
Note that $\la\ge\mu$ implies in particular that $||\la||\ge||\mu||$.
The definition of dominance can be extended to finite non-increasing
sequences of non-negative integers. Observe that adding some zeroes at the end of
two such sequences  $\la$ and $\mu$ does not affect their dominance
relation. 

\begin{Definition}
For an $\ell$-tuple $\la=(\la_1, \ldots,  \la_\ell)$ of non-negative
integers and $k\in \ns$, let $\la^{(k)}$ be the sequence of length $\ell$
obtained by reordering the values $|\la_i-k|$ in non-increasing
order. 
\end{Definition}
\noindent For instance, if $\la=(4,2,1,0)$, then 
$$
\la^{(1)}=(3,1,1,0), \quad 
\la^{(2)}=(2,2,1,0),  \quad   \la^{(3)}=(3,2,1,1), \quad   \la^{(4)}=(4,3,2,0).
$$
The following result was observed experimentally by Jean Creignou in
his work with  Herv\'e Diet on codes in unitary groups and Schur polynomials. More
precisely, it was useful when  implementing their technique to obtain
certain bounds  on the size of codes~\cite{jean}. 

\begin{Proposition}
  Let  $\la=(\la_1, \ldots,  \la_\ell)$ and $\mu=(\mu_1, \ldots,
  \mu_\ell)$ be two partitions of the same weight, possibly completed
  with zeroes so that they have the same length.  Assume that $\la$
  dominates $\mu$. Then for all positive integer $k$,  $\la^{(k)}$
  dominates $\mu^{(k)}$.
\end{Proposition} 

\noindent{\bf Example.} Take  $\la=(4,2,1,0)$ as above and
$\mu=(4,1,1,1)$. These two partitions have the same weight, and
$\la\ge\mu$. The sequences $ \la^{(k)}$ are listed above for $1\le
k\le 4$, and 
$$
\mu^{(1)}= (3,0,0,0), \quad \mu^{(2)}=(2,1,1,1) ,
\quad \mu^{(3)}=(2,2,2,1) ,\quad \mu^{(4)}= (3,3,3,0).
$$
It is easy to check that the proposition holds on this example.

\medskip

Observe that the proposition would be obvious for \emm negative, integers
$k=-m$: in the construction of $ \la^{(k)}$, every part of $\la$  is
simply increased by $m$, and the order of the parts
does not change. For similar reasons, the proposition is clear when
$k\ge \la_1$: the values $\la_i-k$ and $\mu_i-k$ are non-positive for
all $i$, so that the order of the parts is simply reversed:
$\la^{(k)}=(k-\la_\ell,\ldots, k-\la_2, k-\la_1)$.  The fact
that reversion is involved  also shows that the result cannot hold for
partitions of 
different weights. Take for instance $\la=(2)$, $\mu=(1)$ and $k=2$.

\begin{proof}
Recall that $\la$ \emm covers, $\mu$ if there exist no $\nu$ such that
$\la>\nu>\mu$.   The covering relations for the dominance order on
partitions of the same weight were described by Brylawski~\cite{brylawski}. 
The partition
$\la=(\la_1, \la_2,\ldots)$ covers the partition $\mu=(\mu_1, \mu_2, \ldots)$ iff
 there exists $i<j$ such that 
$$
\la_i=\mu_i+1, \la_j=\mu_j -1, \la_p=\mu_p \hbox{ for } p\not \in\{i,j\},
\quad  \hbox{and either } j=i+1 \hbox{ or } \mu_i=\mu_j.
$$

Let us prove the proposition when $\la$ covers
$\mu$. The general case then follows by transitivity.  Recall that $\la^{(k)}$ is obtained by
reordering the multiset $M^{(k)}_{\la}=\{ |\la_p-k|, 1\le p \le \ell\}$. Denote similarly
 $M^{(k)}_\mu=\{ |\mu_p-k| , 1\le p \le \ell\}$. Then $M^{(k)}_\la$ is obtained from $M^{(k)}_\mu$
 by replacing a copy of $ |\mu_i-k|$ by  $ |\mu_i+1-k|$, and a copy of
 $ |\mu_j-k|$ by  $ |\mu_j-1-k|$. 
We study separately 5 cases, depending on the value of $k$. For each of
them, we describe how  $\la^{(k)}$ is obtained from  $\mu^{(k)}$. From
this description, it should be clear that the  dominance relation is preserved.
\begin{itemize}
\item If $k<\mu_j$, the first occurrence of  $\mu_i-k$ in $\mu^{(k)}$ is
  replaced by  $\mu_i-k+1$, while  the last occurrence of  $\mu_j-k$  is
  replaced by  $\mu_j-k-1$. 
\item If $k=\mu_j<\mu_i$, the first occurrence of  $\mu_i-k$ in $\mu^{(k)}$ is
  replaced by  $\mu_i-k+1$, while  the first occurrence of  $0$  is
  replaced by  $1$. 
\item If $k=\mu_j=\mu_i$, the first two copies of  
$0$
  in $\mu^{(k)}$ are both
  replaced by  $1$.
\item If $\mu_j<k\le \mu_i$, then  $\mu^{(k)}$ contains entries
  $\mu_i-k$ and $k-\mu_j$, in some order.  To obtain $\la^{(k)}$, the
  first occurrence of  $\mu_i-k$ is replaced by  $\mu_i-k+1$, and  the
  first occurrence of  $k-\mu_j$ is replaced by   $k-\mu_j+1$.
\item If $k>\mu_i$ , then the values $k-\mu_j$ occur before the values
  $k-\mu_i$ in  $\mu^{(k)}$. To obtain $\la^{(k)}$, the
  first occurrence of  $k-\mu_j$ is replaced by    $k-\mu_j+1$, and  the
  first occurrence of  $k-\mu_i$ is replaced by   $k-\mu_i-1$.
\end{itemize}
\end{proof}

\end{document}